\documentclass[9pt]{article}
\usepackage{amsfonts,amsmath,amssymb,amsthm,tikz}
\usepackage{url}
\usepackage{fullpage}
\usepackage{graphicx}
\usepackage[utf8]{inputenc}
\usepackage[colorlinks=true,citecolor=blue,urlcolor=blue,linkcolor=blue,bookmarksopen=true,unicode]{hyperref}

\usetikzlibrary{calc}

\theoremstyle{plain}
\newtheorem{theorem}{Theorem}
\newtheorem{lemma}[theorem]{Lemma}
\newtheorem{corollary}[theorem]{Corollary}

\newtheorem{proposition}[theorem]{Proposition}

\theoremstyle{remark}

\newtheorem*{remark}{Remark}

\newcommand*{\abs}[1]{\lvert #1\rvert}                           
\newcommand*{\babs}[1]{\bigl\lvert #1\bigr\rvert}                
\newcommand*{\eqdef}{\stackrel{\text{\tiny{def}}}{=}}            

\DeclareMathOperator{\dens}{dens}                                
\DeclareMathOperator{\mindens}{min-dens}                         

\definecolor{grey}{gray}{0.6}

\newcommand\R{\ensuremath{\mathbb{R}}}                           
\newcommand\s{\ensuremath{\mathcal{F}}}                          
\newcommand\Co{\ensuremath{\mathcal{C}}}                         

\newcommand{\pth}[1]{\left( #1 \right)}

\newcommand{\pbrcx}[1]{\ensuremath{\left[ {#1} \right]}}
\DeclareMathOperator{\Expected}{\mathbb{E}}
\newcommand{\Ex}[1]{\Expected\pbrcx{#1}}

\title{Shatter functions with polynomial growth rates}
\author{Boris Bukh\thanks{Department of Mathematical Sciences, Carnegie Mellon University, Pittsburgh, PA 15213, USA. Supported in part by Sloan Research Fellowship and by U.S.\ taxpayers through NSF grant DMS-1301548 and through NSF CAREER grant DMS-1555149.} \and Xavier Goaoc\thanks{Universit\'e Paris-Est, LIGM (UMR 8049), CNRS, ENPC, ESIEE, UPEM, F-77454, Marne-la-Vall\'ee, France. Supported by Institut Universitaire de France.}}
\date{\today}

\begin{document}
\maketitle

\begin{abstract}
  We study how a single value of the shatter function of a set system
  restricts its asymptotic growth. Along the way, we refute a
  conjecture of Bondy and Hajnal which generalizes Sauer's Lemma.
\end{abstract}

\section{Introduction}

A standard tool in combinatorial and
computational geometry is the \emph{shatter function} $f_\s$
of a (geometric) set system $\s$. By \emph{set system} we
mean a family of subsets of a ground set $X$. The \emph{trace} of a
set system $\s$ on a subset $Y\subset X$ is defined as
\[ \s_{|Y} \eqdef \{e \cap Y : e \in \s\}\]
and the \emph{shatter function} of $\s$ is
\[
  f_{\s}(m)\eqdef \max\bigl\{ \abs{\s_{|Y}} : Y\subseteq X, \abs{Y}=m\bigr\}
\]
where $|S|$ denotes the cardinality of a set $S$. The survey of
Matou\v{s}ek~\cite{matousekECM} details several geometric and
algorithmic applications of shatter functions. The
asymptotic growth rate of a shatter function is often its most
important feature.

\bigskip

In this paper, we study how the growth rate of a shatter function can
be controlled by fixing \emph{one} of its values. For example, 
a classical lemma of Sauer~\cite{sauer} (and Vapnik and Chervonenkis~\cite{vc} and
Shelah~\cite{shelah}) asserts that if $f_\s(m)$ is at most $2^m-1$ then $f_\s(n) = O(n^{m-1})$, 
for any natural number $m$. In particular, the growth of a shatter function exhibits
a dichotomy: either $f_{\s}(m)=2^m$ for all $m$, or $f_{\s}$ is bounded by a polynomial. 
We will be concerned below with conditions that ensure precise polynomial growth rates.

\paragraph{New results.}

Let $t_k(m)$ denote the largest integer such that every set system
$\s$ with $f_\s(m) \le t_k(m)$ satisfies $f_\s(n) = O(n^k)$. We prove
the following bounds.

\begin{theorem}\label{t:integer}
  For any integers $m,k \ge 1$,
  \[  (2^{k+1}-k-1)m - 2^{4k} < t_k(m) \le (2^{k+1}-k-1)m+2^{k+1}-k-2.\]
\end{theorem}
We also obtain an analogous result for non-integral values of $k$. The inequalities
are more cumbersome though, see Corollary \ref{cor:irrational} and
Lemma~\ref{l:upper} for the lower and upper bounds respectively.\medskip

We establish the upper bound by a probabilistic construction
(Section~\ref{s:upper}). Interestingly, this upper bound already
refutes a conjecture of Bondy and Hajnal~\cite{bondy}, see
also~\cite[Problem~3.3]{furediP}, regarding a generalization of
Sauer's Lemma: they conjectured that if $f_\s(m) \le g_k(m)$ then $f_\s(n) \le g_k(n)$
for any large enough\footnote{Here ``large enough'' depends solely on
  $m$ and $k$, and not on the set system $\s$. The original statement
  of the conjecture~\cite{bondy} was without this precaution. 
  According to F\"uredi and Pach~\cite[Problem~3.3]{furediP}, the necessity of allowing exceptional values
  is due to Frankl. Bollob\'as and Radcliffe~\cite[Theorem~11]{bollobasR} also           
  gave a probabilistic construction showing that $n_0$ must be larger than $2k$.}~$n$, where
\[ g_k(n) =   1+n+\binom{n}2 + \ldots + \binom{n}{k}.\]
The case $k=m-1$ is Sauer's Lemma and the conjecture was also proven
for $(k,m) = (2,4)$~\cite[Theorem~5]{bollobasR}. If it were true, the
Bondy--Hajnal conjecture would have implied that $t_k(m)$ grows at least as
fast as $\Omega(m^k)$, which is what our upper bound prevents.

We obtain the lower bound by analyzing the density of certain patterns
in simplicial complexes (Section~\ref{s:lower}). This builds on the argument of
Bukh and Conlon~\cite{bukh_conlon} to bound the number of edges of graphs avoiding certain 
subgraphs. Theorem~\ref{t:rational} below specializes to the lower bound in
Theorem~\ref{t:integer} for $s=2^{k+1}-k-1$. 

\begin{theorem}\label{t:rational}
  Let $s \ge 2$ be a rational number, let $q$ be its denominator and
  let $t=\lfloor \log_2 s \rfloor$. Let $m \ge 1$ be an integer. For
  any set system $\s$,
  \[ f_\s(m) \leq s m-3q s^2\log_2 s \quad \implies \quad \forall n \ge m, \quad f_\s(n) \leq
  2^{t+2}m^{2t+2} n^{t+1-(2^{t+1} - t - 2)/(s-1)}.\]
\end{theorem}

As real numbers can be approximated arbitrarily well by rational numbers, it is easy to extend
the preceding theorem to irrational $s$.

\begin{corollary}\label{cor:irrational}
Let $s\geq 2$ be a real number, and let $t=\lfloor \log_2
s\rfloor$. Let $m\geq s^3$ be an integer.  For any set system~$\s$,
  \[ f_\s(m) < s m-10\sqrt{m}\cdot s\sqrt{\log_2 s} \quad \implies \quad \forall n \ge m, \quad f_\s(n) \leq
  3m^{2t} n^{t+1-(2^{t+1} - t - 2)/(s-1)}.\]
\end{corollary}
\begin{proof}
  Let $q=\bigl\lceil(1/s)\sqrt{m/\log_2 s}\bigr\rceil$. Let $s'$ be
  the largest rational number of denominator $q$ such that $s'\leq s$.
  Since $q\leq (2/s)\sqrt{m/\log_2 s}$, we have $3q(s')^2\log_2 s'
  \leq (6/s)\sqrt{m/\log_2 s}\cdot s^2\log_2 s =6\sqrt{m}\cdot
  s\sqrt{\log_2 s}$. Also since $s<s'+1/q$, we have $sm\leq
  s'm+\sqrt{m}\cdot s\sqrt{\log_2 s}$.  Hence $s m-10\sqrt{m}\cdot
  s\sqrt{\log_2 s}<s'm-3q(s')^2\log_2 s'$.
\end{proof}

\paragraph{Related work.}

The only previous lower bound is due to Cheong et
al.~\cite[Theorem~1]{cheongGN}, who proved that $t_k(m) \ge 2^k m -
(k-1) 2^k - 1$ by adapting the inductive proof of Sauer's lemma. The only
upper bound that we are aware of is the easy $t_k(m)=O(m^k/k^k)$, for instance we may
split $[n]$ into $k$ almost equal parts, and let $\s$ consist of those
$k$-sets that contain one vertex from each part.

The argument of Bukh and Conlon was extended from graphs
to hypergraphs by Fitch~\cite{fitch}. Both his and our work use 
generalization of \emph{balanced rooted trees} from the work of 
Bukh and Conlon. There are technical differences, though. 
Fitch works with uniform hypergraphs, whereas we work with simplicial
complexes, which results in a slightly different notion of density.
Our construction (Proposition~\ref{lem:trees}) uses
a different idea from his in \cite[Lemma 1]{fitch}.

\section{Random simplicial complexes}
\label{s:upper}

Recall that a \emph{simplicial complex} is a set system closed under
taking subsets.\footnote{This is usually called an \emph{abstract}
  simplicial complex but since we consider no \emph{embedded}
  simplicial complex in this paper, we do not feel the need to
  emphasize the distinction.} In Lemma~\ref{l:upper} we present a construction of
random simplicial complexes that implies
the upper bound of Theorem~\ref{t:integer}. 

\begin{lemma}\label{l:upper}
  For any real number $s\geq 2$ and for each integer $m\geq 1$, for
  $n$ arbitrarily large there exists a set system $\s$ on $n$ vertices,
  with $f_\s(m) \le sm+(s-1)$ and $f_\s(n) = \Omega\pth{
    n^{t+1-(2^{t+1}-t-2)/(s-1)}}$, where $t=\lfloor \log_2 s\rfloor$.
\end{lemma}
\begin{proof}
  Fix $m$.  For any $n$ large enough, we build a random
  $t$-dimensional simplicial complex $\Co_n$ on $n$ vertices by
  examining each subset of up to $t+1$ vertices in the order of increasing size.
  For each subset $I$, if all $J \subsetneq I$ form faces of our
  complex, we turn $I$ into a face with probability $p=n^{-1/(s-1)}$.
  All choices are independent, and the complex is initialized with all
  $n$ vertices.

  Adding a $t$-dimensional face to $\Co_n$ requires to add each of its
  $2^{t+1} - t - 3$ proper faces of dimension $1$ or more, plus the
  $t$-face itself. The expected number of faces of dimension $t$ of
  $\Co_n$ is thus
  \[ \binom{n}{t+1} p^{2^{t+1}-t-2}.\]
  Let $g(n)$ denote the expected number of faces of $\Co_n$. Note that since $p=n^{-1/(s-1)}$, we have
  $g(n) = \Omega\pth{ n^{t+1-(2^{t+1}-t-2)/(s-1)}}$.

  Let $z=(s-1)(m+1)$. Call an $m$-element set ``bad'' if the set
  contains at least $z$ faces of dimension $1$ or
    more. Since there are at most $2^{2^m}$ complexes on any given
  set of $m$ vertices, the expected number of bad $m$-sets is at most
\[
  2^{2^m}\binom{n}{m} p^z=O(n^m n^{-z/(s-1)})=O(1/n).
\]
  Let $\Co'$ be the complex obtained from $\Co_n$ by removing vertices
  of all bad $m$-sets. Accounting for the traces of size $0$ and $1$, we
    have
\[ f_{\Co'}(m) < z +m+1 = sm+s.\]
  As each vertex belongs to at most $n^t$ faces of $\Co_n$, the
  expected number of faces in $\Co'$ is at least
\[
  g(n)-O(n^{t-1})\geq \tfrac{1}{2}g(n).
\]
So there exists a complex on at most $n$ vertices with at least $\tfrac{1}{2}g(n)$ faces
and $f_{\Co'}(m)\leq sm+(s-1)$. We can ensure that the complex has exactly $n$ vertices by adding
dummy vertices if necessary.
\end{proof}

For any $\varepsilon>0$, Lemma~\ref{l:upper} with
  $s=2^{k+1}-k-1+\varepsilon$ shows that
\[ t_k(m) \le (2^{k+1}-k-1+\varepsilon) m + 2^{k+1}-k-2+\varepsilon.\]
Taking $\varepsilon\to 0$, the upper bound of Theorem~\ref{t:integer} follows.

\begin{remark}
  For most geometric set systems the bound in Sauer's lemma is not sharp. This includes
  the family of halfspaces in $\R^d$. In fact, for this family no shatter
  condition implies the correct bound.

  To see this we may make probabilistic construction
  similar to that of Lemma~\ref{l:upper}. Start with the complete $(d-1)$-dimensional
  skeleton of the $(n-1)$-dimensional simplex, add every $d$-simplex
  randomly and independently, each with probability $p$ with $p =
  \omega(1/n)$ and $p = o(\log n/n)$, then delete every $d$-simplex
  supported on a $m$-element subset of vertices that spans at least
  $m-d+1$ $d$-simplices. With positive probability, the resulting
  random simplicial complex $\Co$ satisfies
  \[f_{\Co}(m) \le 1+m+ \ldots + \binom{m}d + m-d \quad \hbox{and} \quad
   \Ex{f_{\Co}(n)} = \omega(n^d).\]
  Considering points on the moment curve shows that the set system of
  halfspaces in $\R^d$ violates this shatter condition for every $m$.
\end{remark}

\section{Proof of Theorem~\ref{t:rational}}
\label{s:lower}

We first remark that in proving upper bounds on $f_\s(n)$ 
we may restrict ourselves to simplicial
complexes, since any set system can be ``compressed'' without changing
its number of sets nor increasing its shatter function.

\begin{lemma}[Alon~\cite{alon} and Frankl~\cite{frankl}]
  For any finite set system $\s$ there exists an abstract simplicial
  complex~$\Co$ with $|\Co| = |\s|$ and $f_{\Co} \le f_\s$.
\end{lemma}

We write $V(\Co)$ for the set of vertices of a simplicial
complex~$\Co$. Two simplices $\sigma,\sigma'\in \Co$ are
\emph{nonadjacent} in~$\Co$, if they are vertex disjoint, and there is
no edge intersecting both $\sigma$ and $\sigma'$. A set of pairwise nonadjacent
vertices is called an \emph{independent set}. For a complex $\Co$, the
\emph{degree} of a $(d-1)$-simplex $\sigma\in \Co$ is the number of
$d$-simplices $\sigma$ is contained in.  We denote by $\delta_d(\Co)$
the minimum degree of any $(d-1)$-simplex in $\Co$. We define the
  \emph{density} of a subset $S \subseteq V(\Co)$ of vertices of a
  simplicial complex $\Co$ to be $\dens_\Co(S) = e(S)/\abs{S}$, where
  $e(S)$ is the number of non-empty simplices in $\Co$ with at least
  one vertex in~$S$.

\subsection{Balanced rooted $d$-trees and shatter functions}
\label{s:rooted}

A \emph{$d$-tree} is defined inductively. First, a $d$-simplex is a
$d$-tree. If $T$ is a $d$-tree, and $\sigma\in T$ is a
$(d-1)$-simplex, then the complex $T'$ obtained by gluing to $T$
a $d$-simplex formed by $\sigma$ and a new vertex is also a
$d$-tree. We say that $T'$ is obtained by \emph{attaching a vertex to
  $\sigma$}.

\begin{figure}[h]
  \begin{center}  \includegraphics[page=1]{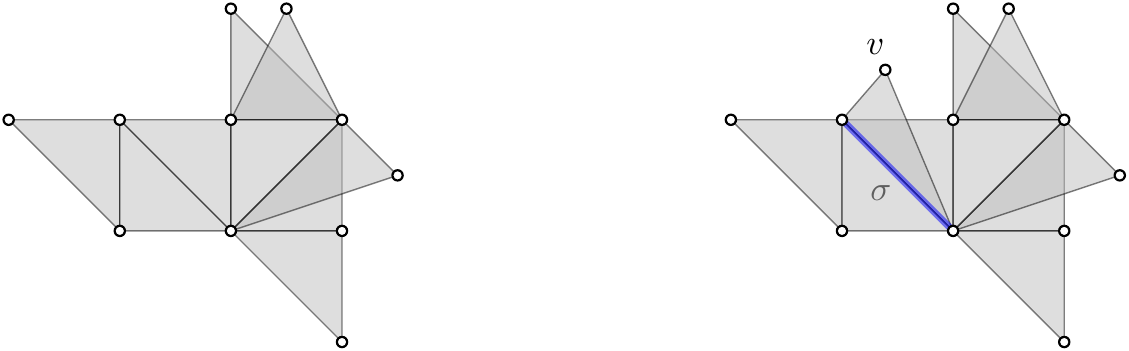} \end{center}
  \caption{Two examples of $2$-trees. The one on the right is obtained by attaching to the one on the left a vertex $v$ to the edge $\sigma$ (in thick blue).}
\end{figure}

A \emph{rooted $d$-tree} $(T,\rho,R)$ consists of a $d$-tree $T$
together with a distinguished $(d-1)$-simplex $\rho\in T$ and an
independent set $R\subset V(T)$ such that $\rho$ is nonadjacent to
each of the vertices in $R$. We call $R$ \emph{vertex roots} and
$\rho$ the \emph{simplex root} of~$T$. The rest of the vertices
we call \emph{unrooted}.

\begin{center}
\begin{tikzpicture}
\filldraw[fill=lightgray] (4.5,1.866025) -- (4,1) -- (5,1) -- cycle;
\filldraw (4.5,1.866025) circle (2pt);
\filldraw[fill=lightgray] (6.866025,0.5) -- (6,1) -- (6,0) -- cycle;
\filldraw (6.866025,0.5) circle (2pt);
\foreach \x in {6,5,...,1}
{
 \filldraw[fill=lightgray](\x,1) -- +(-1,-1) -- +(0,-1) -- cycle;
 \filldraw[fill=lightgray](\x,1) -- +(-1,0) -- +(-1,-1) -- cycle;
 \draw(\x,1) -- +(-1,0);
 \pgfmathtruncatemacro{\lbl}{2*\x+2};
 \filldraw[fill=white] (\x,0) circle (2pt);
 \pgfmathtruncatemacro{\lbl}{2*\x+1};
 \filldraw[fill=white] (\x,1) circle (2pt);
}
\draw[ultra thick] (0,0) -- (0,1);
\filldraw[fill=white] (0,0) circle (2pt); 
\filldraw[fill=white] (0,1) circle (2pt); 
\filldraw (4,0) circle (2pt);
\node at (3,-0.5) {A rooted $2$-tree with three root vertices};
\end{tikzpicture}
\end{center}

The \emph{min-density} $\mindens(T)$ of a rooted $d$-tree $(T,\rho,R)$
is the minimum of $\dens_{T}(S)$ over all non-empty sets $S\subseteq
V(T)\setminus(\rho\cup R)$ of unrooted vertices. If the minimum is
attained by $V(T)\setminus(\rho\cup R)$, the set of all unrooted
vertices, then we call the tree \emph{balanced}.  We use balanced
  rooted trees to bound from below the shatter function of simplicial
  complexes as follows.

\begin{lemma}\label{lem:embed}
  Let $d,m,f,r \ge 1$ be integers. Suppose $(T,\rho,R)$ is a balanced
  $d$-tree with $f$ facets and $r$ vertex roots. Every simplicial
  complex~$\Co$ on $n$ vertices with $\delta_d(\Co)\geq 2m^{1+d/f}
  n^{r/f}$ contains $m$ vertices that span at least
  $\mindens(T)(m-f-d)$ simplices.
\end{lemma}
\begin{proof}
  Let $(T,\rho,R)$ be the $d$-tree in question, and let $\Co$ be a
  simplicial complex on $n$ vertices with $\delta_d(\Co)\geq m^{1+d/f}
  n^{r/f}$.  We assume that $m\geq f+d$, for otherwise the result is
  trivially true, and argue that some $m$-element set $V\subset
  V(\Co)$ spans at least $\mindens(T)(m-f-d)+2^d+r-1$ simplices.

  Fix an arbitrary $(d-1)$-simplex $\sigma$ of $\Co$. Consider copies
  of $T$ such that $\rho$ is mapped to $\sigma$ and different
  $d$-simplices of $T$ are mapped to different $d$-simplices of
  $\Co$. Each such copy can be obtained by embedding facets of $T$
  one-by-one starting with the facet containing the root. Since $T$
  has $f$ facets, then there are at least
  $\delta_d(\Co)(\delta_d(\Co)-1) \ldots (\delta_d(\Co)-f+1) \ge
  (\delta_d(\Co)-f)^f$ copies of the tree $T$ such that $\rho$ is
  mapped to $\sigma$. Since $m\geq f$, it follows that
  $(\delta_d(\Co)-f)^f\geq m^{f+d}n^r$. The pigeonhole principle
  ensures that some $\ell\geq m^{f+d} $ of these copies have the same
  vertex roots; denote them by $T_1,\dotsc,T_{\ell}$.

  Let $V_i=V(T_1)\cup\dotsb\cup V(T_i)$. As $T$ is balanced, $T_1$ has
  at least $\mindens(T)(\abs{V_1}-r-d) + 2^{d}-1 + r$ simplices, and
  at least $\mindens(T)\abs{V_{i+1}\setminus V_i}$ of the simplices
  spanned by $V_{i+1}$ use a vertex from $V_{i+1}\setminus V_i$. Thus,
  by induction on $i$, each $V_i$ spans at least
  $\mindens(T)(\abs{V_i}-r-d)+2^{d}-1 + r$ simplices in~$\Co$.

  The set $V_\ell$ contains at least $\ell$ copies of $T$ with
  prescribed roots, so $m^{d+f}\le \ell \le \binom{\abs{V_\ell}}{\abs{V(T)}} =
  \binom{\abs{V_\ell}}{d+f}$ and it follows that
  $\abs{V_\ell} \ge m$. Since $\abs{V_{i+1}\setminus V_i}\leq f-r$,
  there is $j$ such that $m\geq \abs{V_j}\geq m-(f-r)$. Setting $V=V_j$ we
  obtain the result.
\end{proof}

A similar argument permits us to control the overlap of $d$-simplices.

\begin{lemma}\label{l:overlap}
  Suppose $\Co$ is a simplicial complex, and $\rho$ is a $d'$-simplex in 
  $\Co$ which is contained in $N$ simplices of dimension $d$. Then $\Co$ 
  contains $m$ vertices that span at least 
  $\min\bigl(N,\frac{2^{d+1}-2^{d'+1}}{d-d'}(m-d)\bigr)$ 
  simplices.
\end{lemma}
\begin{proof}
If these $N$ simplices are contained in some $m$-element set, then we are obviously done.
Otherwise, we can find $d$-simplices $\sigma_1,\dotsc,\sigma_{\ell}$ whose union is of
size between $m-d$ and $m$. Let $V_i=\sigma_1\cup\dotsb\cup \sigma_i$.

The number of simplices contained in $\sigma_{i+1}$ that are not contained in $V_i$
is 
\[
  2^{d+1}-2^{\abs{V_i\cap \sigma_{i+1}}}\geq \abs{V_{i+1}\setminus V_i}\frac{2^{d+1}-2^{d'+1}}{d-d'}.
\]
By induction on $i$, it follows that the number of simplices spanned by $V_i$
is at least $\abs{V_i}\frac{2^{d+1}-2^{d'+1}}{d-d'}$.
\end{proof}

\subsection{Construction of balanced $d$-trees of prescribed rational density}

We now prove that balanced $d$-trees of every rational density
exceeding~$2^d$ exist (Proposition~\ref{lem:trees}). The case $d=1$
was previously handled in \cite{bukh_conlon}, and the following
construction borrows some ideas from there. 

\bigskip

Our construction starts with a simplicial complex $T_0$ on the vertex
set $[d(Q+1)]$, whose facets are $d$-simplices of the form
$\{i,i+1,\dotsc,i+d\}$ for all $i=1,\dotsc,dQ$. Alternatively, we can
describe $T_0$ as the complex consisting of all the sets $\sigma
\subset [d(Q+1)]$ satisfying $\max \sigma-\min \sigma\leq d$. For
$i=0,1,\dotsc,Q$ we denote by $\sigma_i$ the $(d-1)$-simplex of $T_0$
defined by
\[ \sigma_i\eqdef\bigl\{ id+1,id+2,\dotsc,(i+1)d\bigr\}.\]
Observe that $(T_0,\sigma_0,\emptyset)$ is a rooted $d$-tree, and that
$\sigma_1,\dotsc,\sigma_Q$ form a partition of unrooted vertices of
this tree. With a slight abuse of notation, we denote this rooted
$d$-tree also by~$T_0$.

\begin{figure}[h]
  \begin{center}  \includegraphics[page=2]{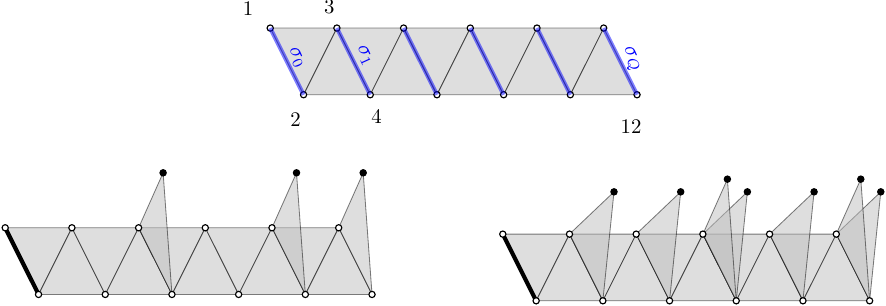} \end{center}
  \caption{The simplicial complex $T_0$ (top) and the rooted $2$-trees $T_3$ (bottom-left) and $T_7$ (bottom-right) for $d=2$ and $Q=5$.}
\end{figure}

If $r<Q$, we define $T_r$ to be the rooted $d$-tree obtained by attaching to $T_0$
a rooted vertex to each of the following $r$ $(d-1)$-simplices
\[ 
  \sigma_{\lceil Q/r\rceil},\sigma_{\lceil 2Q/r\rceil},\dotsc,\sigma_Q.
\]
If $r\geq Q$, we define $T_r$ recursively to be the rooted $d$-tree
obtained from $T_{r-Q}$ by attaching rooted vertices to each of
$\sigma_1,\dotsc,\sigma_Q$.

\begin{proposition}\label{lem:trees}
  For every choice of integers $d,Q\geq 1$, and $r\geq 0$, $T_r$ is a
  balanced $d$-tree with $dQ+r$ facets and $r$ rooted vertices of
  min-density $2^d+\frac{r}{dQ}(2^d-1)$.
\end{proposition}

\noindent
Before we prove Proposition~\ref{lem:trees}, we first argue that the
min-density of $T_r$ is attained on particularly nice sets of unrooted
vertices.

\begin{lemma}\label{lem:nice}
  There exist $1 \le i \le j \le Q$ such that
  $\mindens(T_r)=\dens_{T_r}(\sigma_i \cup \sigma_{i+1} \cup \ldots \cup
  \sigma_j)$.
\end{lemma}
\begin{proof}
  For a set $U$ of unrooted vertices, let the \emph{neighborhood} of $U$ be the set of
  simplices of $T_r$ that contain at least one vertex from~$U$. We denote it by $N(U)$. In
  particular, $e(U)=\abs{N(U)}$. Let $S$ denote a set of unrooted
  vertices that minimizes $\dens_{T_r}(S)$ and is of maximum size among
  such sets.

  \bigskip

  We first claim that $S$ is of the form $S = \sigma_{i_1} \cup \sigma_{i_2}
  \cup \ldots \cup \sigma_{i_p}$. Suppose, for the sake of
  contradiction, that $0<\abs{\sigma_i\cap S}<d$ for some $i$. Pick
  $a,b\in \sigma_i$ such that $a\in S$, $b\notin S$ and $\abs{a-b}=1$.
  By the optimality conditions on $S$, it follows that
  \begin{equation}\label{eq:opt}
    \dens_{T_r}(S)\leq \dens_{T_r}(S\setminus\{a\}), \quad \hbox{and} \quad \dens_{T_r}(S)<\dens_{T_r}(S\cup\{b\}),\end{equation}
  which is equivalent to
  \[ \babs{N(S)\setminus N(S\setminus\{a\})} \leq
  \dens_{T_r}(S), \quad \hbox{and} \quad
  \babs{N(S\cup\{b\})\setminus N(S)}>\dens_{T_r}(S)
  \]
respectively.  It follows that
  $\babs{N(S)\setminus
    N(S\setminus\{a\})}<\babs{N(S\cup\{b\})\setminus N(S)}$. To reach
  a contradiction we now exhibit an injective map from
  $N(S\cup\{b\})\setminus N(S)$ to $N(S)\setminus N(S\setminus\{a\})$.

  Assume that $b=a+1$, for the other case $b=a-1$ is analogous, and
  consider the map
  \[
  \phi(\sigma)\eqdef
  \begin{cases}
    \bigl(\sigma\cup\{a\}\bigr)\setminus \{a+1\} & \text{if }a+d+1\notin \sigma,\\
    \bigl(\sigma\cup\{a\}\bigr)\setminus \{a+d+1\} & \text{if }a+d+1\in
    \sigma.
  \end{cases}
  \]
  If $\sigma \in T_r$, that is $\max \sigma-\min \sigma\leq d$, and $b
  = a+1 \in \sigma$ then $\max \sigma \le a+d+1$ and it follows that
  $\phi(\sigma) \in T_r$. Moreover, if $\sigma \cap S =\emptyset$ then
  $\phi(\sigma) \cap S = \{a\}$. This implies that $\phi$ maps $
  N(S\cup\{b\})\setminus N(S)$ to $N(S)\setminus N(S\setminus\{a\})$;
  this map is easily seen to be injective. The existence of $\phi$
  contradicts the optimality of $S$, and thus each $\abs{S \cap
    \sigma_i}$ is $0$ or $d$.
  
  \bigskip

  We can now partition $S = S_1 \cup S_2 \cup \ldots \cup S_{p'}$
  where each $S_i$ is a maximal union of consecutive
  $\sigma_j$'s. Since
  \[
  e(S)=e(S_1\cup\dotsb\cup S_{p'})=e(S_1)+\dotsb+e(S_{p'})\geq
  \bigl(\abs{S_1}+\dotsb+\abs{S_{p'}} \bigr)\min_i
  \dens_{T_r}(S_i)=\abs{S}\min_i\dens_{T_r} (S_i) ,
  \]
  it follows that $\min_i \dens_{T_r}(S_i)\leq
  \dens_{T_r}(S)=\mindens(T_r)$. This completes the proof.
\end{proof}

\begin{proof}[Proof of Proposition~\ref{lem:trees}]
  In view of Lemma~\ref{lem:nice}, it remains to compute
  $\dens_{T_r}(\sigma_i \cup \sigma_{i+1} \cup \ldots \cup \sigma_j)$
  and show that it is minimal for $i=1$ and $j=Q$. Let $S=\sigma_i\cup
  \sigma_{i+1}\cup \dotsb\cup\sigma_j$. For $r\ge Q$,
  \begin{equation}\label{eq:propagate}
    \dens_{T_r} (S)=\dens_{T_{r-Q}}(S)+(2^d-1)/d,
  \end{equation}
  so we focus on the cases $r<Q$.

  \bigskip

  We first express $\dens_{T_r} (S)$ in terms of $L[i,j]$, where for
  $1 \le i \le j \le Q$
  \[ L[i,j]\eqdef \babs{\{i,i+1,\dotsc,j\}\cap\{\lceil
    Q/r\rceil,\lceil 2Q/r\rceil,\dotsc,Q\}}. \]
  The computations are easiest when $j=Q$. In this case, whenever a
  simplex $\sigma\in T_0$ meets $S$, we necessarily have $\max \sigma
  \in S$.  For each $x\in S$ there are exactly $2^d$ simplices
  $\sigma\in T_0$ such that $\max \sigma=x$.  Therefore the number of
  simplices of $T_0$ that meet $S$ is $2^d\abs{S}$.  Adding the
  simplices contained in the facets counted by $L[i,j]$ we obtain
  \[
  \dens_{T_r}(S)=\frac{2^d\abs{S}+(2^d-1)L[i,j]}{\abs{S}}=2^d+\frac{(2^d-1)L[i,Q]}{d(Q-i+1)}\qquad\text{for
  }S=\sigma_i\cup\sigma_{i+1}\cup\dotsb\cup \sigma_Q.
  \]

  When $j<Q$, the computation is similar except that we also need to
  count the simplices $\sigma$ of $T_0$ such that $\max \sigma\notin
  S$. Call such a simplex $\sigma$ \emph{dangling}.  If $\sigma$ is
  dangling, then $\min \sigma\in \sigma_j$.  Furthermore, for each
  $\ell\in\{1,2,\dotsc,d\}$ there are $2^d$ simplices $\sigma$
  such that $\min \sigma=jd+\ell$ and exactly $2^d-2^{d-\ell}$
  of them are dangling. The
  total number of dangling simplices is thus
  \[
  \sum_{\ell=1}^d \bigl(2^d-2^{d-\ell}\bigr)=d2^d-(2^d-1)=(d-1)2^d+1,
  \]
  yielding
  \[
  \dens_{T_r}(S)=2^d+\frac{(2^d-1)L[i,j]+(d-1)2^d+1}{d(j-i+1)}\qquad
  \text{for }S=\sigma_i\cup \sigma_{i+1}\cup \dotsb\cup \sigma_j,\
  j<Q.
  \]

  \bigskip

  Next, note that if $i>1$ and $L[i,j]=L[i-1,j]$, then $\dens_{T_r}(S\cup
  \sigma_{i-1})<\dens_{T_r}(S)$. Similarly, if $j\neq Q$, and
  $L[i,j]=L[i,j+1]$, then $\dens_{T_r}(S\cup \sigma_{j+1})<\dens_{T_r}(S)$. As we
  look for the minimum density, we may assume that
  \[ i\in \{1, \lceil Q/r\rceil+1,\lceil 2Q/r\rceil+1,\dotsc,\lceil
  (r-1)Q/r\rceil+1\}, \quad \hbox{and} \quad j\in \{\lceil
  Q/r\rceil-1,\lceil 2Q/r\rceil-1,\dotsc,Q-1,Q\}.\]
  We can thus assume that $i=\lceil a Q/r\rceil+1$ with $0\leq a <r$
  and that either $j=Q$ or $j = \lceil bQ/r\rceil-1$ with $a < b \le
  r$.

  \bigskip

  If $j=Q$ then
  \begin{align*}
    \dens_{T_r}(S)=2^d+\frac{(2^d-1)(r-a)}{d(Q-\lceil a
      Q/r\rceil)}=2^d+\frac{(2^d-1)(r-a)}{d(\lfloor
      (r-a)Q/r\rfloor)}\geq 2^d+\frac{r(2^d-1)}{dQ},
  \end{align*}
  with equality if $a=0$. If $j = \lceil bQ/r\rceil-1$ then
\begin{align*}
  \dens_{T_r}(S)&=2^d+\frac{(2^d-1)(b-a-1)+(d-1)2^d+1}{d(j-i+1)}\\
          &=2^d+\frac{(2^d-1)(b-a)+(d-2)2^d+2}{d(j-i+1)}\\
          &>2^d+\frac{(2^d-1)(b-a)}{d(bQ/r -  a Q/r)}+\frac{(d-2)2^d+2}{d(j-i+1)}\\
          &=2^d+\frac{2^d-1}{dQ/r}+\frac{(d-2)2^d+2}{d(j-i+1)}
\end{align*}
As this exceeds $2^d+\frac{2^d-1}{dQ/r}$, we conclude that
the min-density of $T_r$ is attained for $S= \sigma_1 \cup \sigma_2
\cup \ldots \cup \sigma_Q$ and is
$\mindens(T_r)=2^d+\frac{2^d-1}{dQ/r}$.
\end{proof}

\subsection{Wrapping up} 

We can now prove Theorem~\ref{t:rational}. Let $s \ge 2$ be a rational
number, let $q$ denote its denominator, and let $t=\lfloor \log_2 s
\rfloor$. We want to prove that for any simplicial complex $\Co$,
\[ f_\s(m) \le s m-3qs^2 \log_2 s \quad \Rightarrow \quad f_\s(n) \le 
2^{t+2}m^{2t+2} n^{t+1-(2^{t+1} - t - 2)/(s-1)}.\]
Assume that $\Co$ has more than $2^{t+2}m^{2t+2} n^{t+1-(2^{t+1} - t -
  2)/(s-1)}$ simplices on $n$ vertices.  We will show that some $m$ of
its vertices must span more than $sm-3q s^2\log_2 s$ simplices. We assume that $m \ge 3s\log_2 s$ as otherwise the statement
  holds trivially. 

For each $0\leq d\leq t$, let $t_d=(s-2^d)/(s-1)$. Let
$s_d=1+t_1+t_2+\dotsb+t_d$ and note that $s_d= d+1 - (2^{d+1} - d -
2)/(s-1)$. In particular, our assumption states that $\Co$ contains
more than $2^{t+2}m^{2t+2} n^{s_t}$ simplices. We note that $s_t
\ge t+1-\frac{2s-t-2}{s-1} > t-1$.

\paragraph{Case 1.}

We first consider the case where there is some $d$ with $1 \le d \le
t$ such that $\Co$ contains more than $2^dm^{2d}n^{s_d}$ simplices of
dimension $d$. Pick the smallest such $d$.

Now, delete from $\Co$ all $(d-1)$-simplices of degree less than
$2m^2n^{t_d}$, also removing any simplices that contain them. By
minimality of $d$, doing so removes fewer than
$2^{d-1}m^{2d-2}n^{s_{d-1}}\cdot 2m^2 n^{t_d}$ of the $d$-simplices.
Hence, the resulting complex $\Co'$ contains at least one $d$-simplex,
and satisfies $\delta_d(\Co')\geq 2m^2 n^{t_d}$.

Write the rational number $(t_d^{-1}-1)/d=(2^d-1)/d(s-2^d)$ in the
form $Q/r$ with $\gcd(Q,r)=1$. Let $T$ be a balanced rooted $d$-tree
with parameters $Q$ and $r$ as given by Proposition~\ref{lem:trees}.
Note that $\mindens T=s$. By Lemma~\ref{lem:embed} there is set of $m$
vertices on $\Co'$ (and hence of $\Co$) that spans at least
$s(m-d(Q+1)-r)$ simplices. The requisite bound follows from $Q\leq
2^dq\leq sq$ and $r\leq d s q$, and from $d\leq t\leq \log_2 s$.

\paragraph{Case 2.}

In the remaining case, the number of simplices of dimension up to $t$
is at most
\[ \sum_{d=0}^t 2^dm^{2d}n^{s_d} \leq \sum_{d=0}^t 2^d m^{2t} n^{s_t} \leq 2^{t+1}m^{2t}n^{s_t}.\]
so $\Co$ contains more than $2^{t+1}m^{2t+2}n^{s_t}$ simplices
of dimension greater than $t$.

If $\Co$ contains a $d$-simplex for some $d\geq 2\log_2 m$, then any
$m$-set containing this simplex contains at least $2^d\geq m^2>sm$
simplices.  So assume that $\Co$ is of dimension at most $2\log_2 m$.

By the pigeonhole principle, there is a $d$ with $t+1 \le d \le 2
  \log_2 m -1$ such that $\Co$ contains at least $2^{t+1}
m^{2t+2}n^{s_t}/(2\log_2 m-t-1)\geq
2^{t+1}(m^{2t+2}/2\log_2m) n^{s_t}$ simplices of dimension $d$.

Since $\Co$ contains at most $m^{2t}n^{s_t}$ simplices of dimension
$t$, it follows by another application of the pigeonhole principle
that there is a $t$-simplex $\rho$ that is contained in at least
$2^{t+1}m^2/2\log_2 m$ simplices of dimension $d$.  Lemma~\ref{l:overlap} then
yields an $m$-element set spanning at least
\[
  \min\bigl(2^{t+1}m^2/2\log_2 m,\frac{2^{d+1}-2^{t+1}}{d-t}(m-d)\bigr)
\]
simplices. We have $2^{t+1}m^2/2\log_2 m\geq sm$. Also,
$\frac{2^{d+1}-2^{t+1}}{d-t}(m-d) \geq 2^{t+1}(m-t-1)\geq
sm-s\log_2s$. 


\bibliographystyle{alpha}
\bibliography{ref}
\end{document}